\documentclass{amsart}
\usepackage[utf8]{inputenc}
\usepackage{amsmath}
\usepackage{mathtools}

\theoremstyle{definition}
\newtheorem{theorem}{Theorem}[section]

\newtheorem{definition}[theorem]{Definition}
\newtheorem{remark}[theorem]{Remark}
\newtheorem{proposition}[theorem]{Proposition}
\newtheorem{example}[theorem]{Example}
\newtheorem{corollary}[theorem]{Corollary}

\title{Lipschitz normally embedded sets do not need to have Lipschitz normally embedded medial axis}
\author{Michał Kosiba}
\address{Jagiellonian University, Faculty of Mathematics and Computer Science, Łojasiewicza 6, 30-348 Kraków, Poland}
\email{michal.kosiba@doctoral.uj.edu.pl}

\begin{document}

\begin{abstract}
The aim of this paper is to study the Lipschitz normally embedded property for a set and its medial axis. We consider if and when a non-LNE set implies non-LNE medial axis and converse. We present an example a of Lipschitz normally set that has medial axis which is not Lipschitz normally emebedded. At the end we discuss special case when a set is a one dimensional germ on a plane.
\end{abstract}

\maketitle

\section{Introduction}
Let us consider a path connected semialgebraic set $X \subset{R}^{n}$. Given two points $x,y \in X$ we can calculate distance between them in two natural way. The first approach is to consider the Eulidean distance induced from $\mathbb{R}^{n}$, and the second idea is two calculate infimum of rectifiable paths contained in $X$ connecting $x$ and $y$. In such way we can introduce two metrics on $X$ - induced metric (Euclidean) and inner metric (given by legths of curves). We call set $X$ Lipschitz normally embedded - LNE (as introduced by L. Birbrair and T. Mostowski) if these two metrics are equivalent. The aim of this paper is to answer the question if the information of LNE property of a set is hold by its medial axis.

For a set $X \subset \mathbb{R}^{n}$ its medial axis is a collection of this points for which the distance to $X$ is realized by more than one point in $X$. Although having a simpler structure, medial axis was considered to hold all the most important information of its set. We could expect that non-LNE medial mean that the set itself would not be LNE but the situation turns out to be more subtle. Only on a plane non-LNE medial axis implies non-LNE set. In general such implication is not true and we will present an example of such a set.

\section{Lipschitz normally embedded sets - general informations}

Let $X \subset \mathbb{R}^{n}$ be a connected semialgebraic set. Although the normal embeddings are usually studied for semialgebraic sets, all the results below generalize to sets definable in polynomially bounded o-minimal structures and hence remain true for subanalytic germs.

\begin{definition}
The \textit{induced metric} on $X$, which we will denote by $d_{ind}$, is the one coming from the Euclidean norm, i.e. for $x_{1},x_{2} \in X$, $d_{ind}(x_{1},x_{2})=||x_{1}-x_{2}||$.
\end{definition}

\begin{definition}
Let $x_{1},x_{2} \in X$ and let $\Gamma_{x_{1},x_{2}}=\{\gamma:[0,1] \rightarrow X\;|\; \gamma(0)=x_{1},\; \gamma(1)=x_{2},\; \gamma \;\text{rectifiable}\}$. We define the \textit{length metric} (or {\it inner metric}) on $X$ as $d_{l}(x_{1},x_{2})=\inf_{\gamma \in \Gamma_{x_{1},x_{2}}} l(\gamma)$, where $l(\gamma)$ means the length of $\gamma$.
\end{definition}

\begin{remark}
The definition above is correct since for semialgebraic sets being connected is equivalent to being path connected and this, moreover, by rectifiable paths (it follows from cylindrical cell decompositions).
\end{remark}

Notice that for $x_{1},x_{2} \in X$ we always have $d_{ind}(x_{1},x_{2}) \leq d_{l}(x_{1},x_{2})$. There is a natural question when this can be reversed with a constant. That motivates the following definiton.

\begin{definition}
The set $X$ is called \textit{Lipschitz normally embedded (LNE)} in $\mathbb{R}^{n}$ if the metrics $d_{ind}$ and $d_{l}$ are equivalent, which means that there exists $C>0$ such that for $x_{1},x_{2} \in X$ holds $d_{l}(x_{1},x_{2}) \leq Cd_{ind}(x_{1},x_{2})$. 
\end{definition}

We have a very nice theorem concerning the possibility of decomposing a definable into normally embedded components.

\begin{theorem}[Pancake decomposition, \cite{normal}, \cite{kurdyka}, \cite{kurdmost}, \cite{paru}]
\label{pancake}
Let $X \subset \mathbb{R}^{n}$ be a closed semialgebraic set. Then there exists a finite family of subsets $\{X_{i}\}_{i=1}^{k}$ called {\it pancakes} such that:\\
1) all $X_{i}$ are definable closed subsets of $X$;\\
2) $X=\bigcup_{i=1}^{k} X_{i}$;\\
3) $\text{dim}(X_{i} \cap X_{j}) < \min(\text{dim}X_{i},\text{dim}X_{j})$, for every $i \neq j$;\\
4) all $X_{i}$ are normally embedded in $\mathbb{R}^{n}$.
\end{theorem}

\begin{definition}
Let $X \subset \mathbb{R}^{n}$ be a closed connected semialgebraic set with its pancake decomposition $\{X_{i}\}_{i=1}^{N}$. Take $x_{1},x_{2}$ and consider a sequence of points $\{y_{1},...,y_{k}\}$ that satisfies the following conditions:\\
1) $y_{1}=x_{1}$ and $y_{k}=x_{2}$;\\
2) every pair $y_{j},y_{j+1}$ lies in one pancake $X_{i}$;\\
3) if $y_{j},y_{j+1} \in X_{i}$, then $y_{s} \notin X_{i}$ for $s \notin \{j,j+1\}$.\\
Let $Y_{x_{1},x_{2}}$ be the set of all finite sequences satisfying the conditions above. For every sequence $y=\{y_{1},...,y_{k}\} \in Y_{x_{1},x_{2}}$, we set $l(y)= \sum_{j=2}^{k}d_{ind}(y_{j},y_{j-1})$ and $d_{P}(x_{1},x_{2})=\inf_{y \in Y_{x_{1},x_{2}}} l(y)$.
\end{definition}

\begin{theorem}[\cite{normal} Theorem 3.1]
The function $d_{P}:X \times X \rightarrow \mathbb{R}$ is semialgebraic and defines a metric on $X$.
\end{theorem}

Such a metric as in the theorem above will be called a $\textit{pancake metric}$. The interesting property of this metric is the fact that the distance is always realized by some sequence of points as stated in the proposition below.

\begin{proposition}[\cite{normal} Lemma 3.2]
For $x_{1},x_{2} \in X$, there exists $y \in Y_{x_{1},x_{2}}$ such that $d_{P}(x_{1},x_{2})=l(y)$. 
\end{proposition}

The next important property of the pancake metric is the fact that it is equivalent to the length metric. This means that we can study connected semialgebraic sets using the pancake metric (which is definable) instead of the length metric.

\begin{theorem}[Kurdyka-Orro, see also \cite{normal} Theorem 3.2]
The pancake metric is equivalent to the length metric.
\end{theorem}

\begin{definition}
Let us recall that for a continuous, nonconstant, semialgebraic function germ $f\colon (\mathbb{R}_+,0) \rightarrow (\mathbb{R},0)$ (here $f(0)=0$), we have a representation of the form $f(t)=at^{\alpha}+o(t^{\alpha})$ for some $a\in \mathbb{R}\setminus\{0\}$ and $\alpha \in \mathbb{Q}$. The number $\alpha>0$ is called $\textit{the order of $f$ at 0}$ and denoted $ord\ f$.
\end{definition}

\begin{definition}
Let $\gamma_{i}:[0,\epsilon) \rightarrow \mathbb{R}^{n}$ be semialgebraic curves for $i=1,2$ such that $\gamma_{1}(0)=\gamma_{2}(0)=0$, and $\gamma_{i} \subset X$.\\
We define the \textit{outer order of tangency} as:
\[tord(\gamma_{1},\gamma_{2})=ord\ d_{ind}(\gamma_{1}(t),\gamma_{2}(t)),\]
and the \textit{inner order of tangency} by
\[tord_{inn}(\gamma_{1},\gamma_{2})=ord\ (d_{P}(\gamma_{1}(t),\gamma_{2}(t))).\]
\end{definition}

\begin{remark}
Although the pancake decomposition for a set could not be unique, the order of pancake metric does not depend on the choice of decomposition. It follows from the equivalence with the length metric.
\end{remark}

\begin{definition}
Let $\gamma :[0,\epsilon) \longrightarrow \mathbb{R}^{n}$ such that $\gamma(0)=0$. We say that $\gamma$ is \textit{parametrized by the distance} if $||\gamma(t)||=t$ for every $t \in [0,\epsilon)$.
\end{definition}

\begin{remark}
Notice that every definable curve can be parametrized by the distance. This is essentially due to the fact that a one dimensional definable germ cannot `oscillate' so that for small arguments $t$ the intersection of $\gamma([0,\epsilon))$ with the sphere centred at 0 with radius $t$ is a singleton.
\end{remark}

\begin{theorem}[\cite{arccrit}, Theorem 2.2]
\label{lne}
Let $X$ be a closed semialgebraic germ at $0\in\mathbb{R}^n$. The following assertions are equivalent:\\
i) the germ of $X$ at $0$ is normally embedded;\\
ii) there exists a constant $C>0$ such that for any pair of arcs $\gamma_{1},\gamma_{2}$, parametrized by the distance at $0$, $(\gamma_{i}(0)=0)$ we have:
\[d_{l}(\gamma_{1}(t),\gamma_{2}(t)) \leq C d_{ind}(\gamma_{1}(t),\gamma_{2}(t));\]
iii) for any pair of arcs $\gamma_{1},\gamma_{2}$, parametrized by the distance to $0$ we have:
\[tord(\gamma_{1},\gamma_{2})=tord_{inn}(\gamma_{1},\gamma_{2}).\]
\end{theorem}

Notice that if we have two semialgebraic curves $\gamma_{i}(t):[0,\epsilon] \rightarrow \mathbb{R}$, then the functions $d_{ind}(\gamma_{1}(t),\gamma_{2}(t))$ and $d_{P}(\gamma_{1}(t),\gamma_{2}(t))$ are semialgebraic and obviously if for some $t_{0}$ we have $d_{ind}(\gamma_{1}(t_{0}),\gamma_{2}(t_{0})=0$, then $d_{P}(\gamma_{1}(t_{0}),\gamma_{2}(t_{0}))=0$. Restricting the functions to the images of $(\gamma_{1},\gamma_{2})$, we have $d_{ind}^{-1}(0) \subset d_{P}^{-1}(0)$ which means that the assumptions of the classical version of Łojasiewicz's inequality are satisfied and we have $\alpha,C>0$ such 
\[d_{P}(\gamma_{1}(t),\gamma_{2}(t)) \leq Cd_{ind}(\gamma_{1}(t),\gamma_{2}(t))^\alpha. \;\;\; (**)\]
From the fact that $d_{ind}(\gamma_{1}(t),\gamma_{2}(t)) \leq d_{P}(\gamma_{1}(t),\gamma_{2}(t))$ we have that $\alpha \leq 1$.

\begin{definition}
Let $X$ be a closed semialgebraic set, $0 \in X$ a fixed point and let $\gamma_{i}:[0,1) \rightarrow X, \; i=1,2$ be two semialgebraic curves such that $\gamma_{i}(0)=0$ and giving distinct germs at 0.
We define the Łojasiewicz exponent for a such pair of curves as
\[\mathcal{L}_X(\gamma_{1},\gamma_{2})=\inf\{\alpha \geq 1 \;|\; (**)\; \text{holds for} \; \alpha  \; \text{and some} \; C>0\}\]
which is attained (by the Bochnak-Risler Theorem or, in this case also the Curve Selection Lemma),
and the Łojasiewicz's exponent for the closed semialgebraic germ $(X,0)$  as:
\[\mathcal{L}(X,0)=\sup \{\mathcal{L}_X(\gamma_{1},\gamma_{2})\;|\; \gamma_{i}:[0,\epsilon) \rightarrow X,\; \gamma_{i}(0)=0,\; \gamma_{i} \;\text{is parametrized by the distance}\}.\]
\end{definition}

\begin{theorem}
Let $(X,0)$ be a closed semialgebraic germ. Then:
\[(X,0) \; \text{is normally embedded} \Longleftrightarrow \mathcal{L}(X,0)=1 \]
\begin{proof}
Follows directly from the Theorem \ref{lne}. 
\end{proof}
\end{theorem}

\section{Link Lipschitz normally embedded sets}

Sometimes, studying a LNE property for a set in general can be difficult. However, it turns oout that there is a direct equivalence between LNE property of the set and LNE property of its sections. It means we can determine if the set is LNE by checking the LNE property of its sections.

\begin{definition}
Let $X \subset \mathbb{R}^{n}$ be a subset, $p \in \overline{X}$ and $X_{t}:=X \cap \mathbb{S}^{n-1}(p;t)$ for all $t \geq 0$. We say that $X$ is \textit{link Lipschitz normally embedded (at p) (LLNE)} if there is a constant $C \geq 1$ such that $d_{X_{t}} \leq C ||\cdot||$, for $t>0$ small enough.  
\end{definition}

\begin{definition}
Let $||\cdot||_{1}$ be a subanalytic norm on $\mathbb{R}^{n}$. We say that $X$ is \textit{LLNE with respect to (w.r.t.) $||\cdot||_{1}$} if there is a constant $C>0$ such that $d_{X_{t}^{||\cdot||_{1}}} \leq C||\cdot||$, for all small enough $t>0$, where $X_{t}^{||\cdot||_{1}} = \{x \in X \; |\; ||x-p||_{1}=t\}$.
\end{definition}

\begin{definition}
Given $v = (v_{1}, ..., v_{n}) \in \mathbb{R}_{+}^{n}= \{(x_{1},...,x_{n}); x_{i}>0, \forall i\}$ we define the subanalytic norm $||x||_{max,v}=\{v_{1}|x_{1}|,...,v_{n}|x_{n}|\}$.    
\end{definition}

The next proposition shows that LLNE property does not depend on the choice of a subanalytic norm

\begin{proposition}[\cite{link}]
Let $X$ be a subanalytic set, $p \in \overline{X}$. Let $||\cdot||_{1}$ and $||\cdot||_{2}$ be a subanalytic norms on $\mathbb{R}^{n}$. Then $X$ is LLNE at p w.r.t. $||\cdot||_{1}$ if and only if $X$ is LLNE at p w.r.t $||\cdot||_{2}$.   
\end{proposition}

Finally we can state the equivalence of LNE and LLNE properties.

\begin{theorem}[\cite{link}]
\label{LLNE}
Let $X \subset \mathbb{R}^{N}$ be a closed subanalytic set, $0 \in X$. Let $C_{1},...,C_{r}$ be the connected components of $X \backslash \{0\}$ (as a germ). Then the following are equivalent:
\begin{enumerate}
    \item $X$ is LNE at $0$;
    \item Each $\overline{C_{i}}$ is LNE at $0$ and there exists $K>0$ such that $d_{0}(X_{t}) \geq Kt$ for all small enough $t>0$;
    \item Each $\overline{C_{i}}$ is LLNE at $0$ and there exists $K>0$ such that $d_{0}(X_{t}) \geq Kt$ for all small enough $t>0$.
\end{enumerate}
\end{theorem}

\section{Lipschitz normally embedded sets and medial axis}

\begin{definition}
Let $X\subset \mathbb{R}^{n}$ be a closed set and let $x \in \mathbb{R}^{n}$. We define \textit{the set of the closest points to $x$}:
$$m(x)=\{y\in X \;|\; ||y-x||=d(x,X)\}.$$
Due to the fact that $X$ is closed, $m(x)$ is nonempty for every $x \in \mathbb{R}^{n}$. Moreover there holds
$$m(x)=X \cap \mathbb{S}(x,d(x,X)),$$
\end{definition}

After defining this object we can finally introduce a concept of medial axis.

\begin{definition}
For a closed nonempty set $X \subset \mathbb{R}^{n}$ we define the set
$$M_{X}=\{x \in \mathbb{R}^{x} \backslash X\;|\; \# m(x)>1\}$$
and call it the \textit{medial axis of the set} $X$.
\end{definition}

Medial axis is a simpler object that still holds many informations about the set. There is a natural qurstions if it preserves the LNE property. The answer turns out to be negative as show in examples below. Examples of three cases are rather simple:\\
- Non-LNE set with LNE medial axis: a cusp on a plane;\\
- LNE set with LNE medial axis: a graph of absolute value of real numbers;\\
- Non-LNE set with non-LNE medial axis: a germ of three curves tanget to each other (described in details in Section 5)\\

The case of LNE set with non-LNE medial axis is not as obvious but such an object exists too as shown in example below

\begin{example}
Let us consider four curves in $\mathbb{R}^{3}$: \\
$\gamma_{11}=\{(x,y,0) \in [0,1]\times\mathbb{R}_{+}\times\{0\}\;|\; y=\sqrt{x}\}$, \\
$\gamma_{12}=\{(x,y,0) \in [0,\frac{1}{2}]\times\mathbb{R}_{+}\times\{0\}\;|\; y=2\sqrt{x}\}$, \\
$\gamma_{21}=\{(x,y,) \in [-1,0]\times\mathbb{R}_{+}\times\{0\}\;|\; y=\sqrt{-x}\}$, \\
$\gamma_{21}=\{(x,y0) \in [-\frac{1}{2},0]\times\mathbb{R}_{+}\times\{0\}\;|\; y=2\sqrt{-x}\}$.\\
Now define following three sets:\\
$X_{1}=\{(x,y,z) \in \mathbb{R}\times[0,1]\times \mathbb{R}\;|\; (x_{1},y_{1},0) \in \gamma_{11}, (x_{2},y_{2},0) \in \gamma_{12} \; \text{with} \; y=y_{1}=y_{2} \Longrightarrow ||\frac{y_{1}+y_{2}}{2} -y_{2}||=||\frac{y_{1}+y_{2}}{2} -y||\}$\\
$X_{2}=\{(x,y,z) \in \mathbb{R}\times[0,1]\times \mathbb{R}\;|\; (x_{1},y_{1},0) \in \gamma_{21}, (x_{2},y_{2},0) \in \gamma_{22} \; \text{with} \; y=y_{1}=y_{2} \Longrightarrow ||\frac{y_{1}+y_{2}}{2} -y_{2}||=||\frac{y_{1}+y_{2}}{2} -y||\}$\\
$X_{3}=\{(x,y,z) \in [-\frac{1}{2},\frac{1}{2}] \times [0,1] \times \{0\} \; | |y| >2\sqrt{|x|}\;\}$.\\
It is easy to see that $X_{1},X_{2}$ are just horns created respectively on $\gamma_{11}, \gamma_{12}$ and $\gamma_{21}, \gamma_{22}$, and $X_{3}$ is some kid of wall joining them.\\
Finally we define set $X = X_{1} \cup X_{2} \cup X_{3}$.
\end{example}

\begin{proposition}
Set X in example above is LNE but its medial axis is non-LNE.
\begin{proof}
First notice that medial axis of $X$ is non-LNE. The components of the medial axis contained inside the horns are just curves that has order of tangency greater than one. On the other hand, they meet only at the orgin so the inner distance has order 1, and the medial axis cannot be LNE.\\
To prove that $X$ is LNE let us consider the maximum norm in $\mathbb{R}^{3}$.
Notice that $X_{i}^{\max}=X \cap \mathbb{S}^{2}_{\max}(0, t)$ for $t\in (0,1)$ are just two circles connected with an interval. Such a section is LNE with constant 1. Then from Theorem \ref{LLNE} we obtain that $X$ is LNE set.
\end{proof}
\end{proposition}

\section{Plane case}

There is a natural question if there are some special cases in which non-LNE medial axis implies the lack of the LNE property of the set. In this section we prove that in $\mathbb{R}^{2}$ a non-LNE medial axis (as a germ) implies a non-LNE set (as a germ)

\begin{definition}
Let $X \subset \mathbb{R}^{n}$ be a closed set, $a \in X$. We define the \textit{Peano tangent cone} at $a$ as:
\[C_{a}(X)=\{v \in \mathbb{R}^{n}\;|\; \exists\;X\ni x_{\nu} \rightarrow a,\; \exists\; t_{\nu}>0,\;t_{\nu}(x_{\nu}-a) \rightarrow v \}.\]

\end{definition}
It is known that for a non-constant definable curve germ $\gamma : [0,\varepsilon)\to\mathbb{R}^n$ identified with its image, $C_0(\gamma)$ is a half-line.

\begin{theorem}
\label{curvegerm}
Let $\gamma_{1},\gamma_{2}:[0,1) \rightarrow \mathbb{R}^{n}$ be two semialgebraic curves such that $\gamma_{1}(0)=\gamma_{2}(0)=0$. We identify them with their images. Let $l_{1},l_{2}$ be the tangent half-lines at $0$ respectively to $\gamma_{1}$ and $\gamma_{2}$. Then
\[\text{the germ}\; (\gamma_{1} \cup \gamma_{2},0)\; \text{is normally embedded} \Longleftrightarrow l_{1} \neq l_{2}.\]
\begin{proof}
From Theorem \ref{lne} we know that we have to check whether the outer and inner orders of tangency of all pairs of curves parametrized by the distance are equal. Since we have a germ of two curves, we have to check only one pair of curves -- $(\tilde{\gamma_{1}},\tilde{\gamma_{2}})$ where $\tilde{\gamma_{i}}$ is a reparametrization of $\gamma_{i}$ such that it is parametrized by the distance. We may assume this was the case for the initial curves. Let $v_{i} \in C_{0}(\gamma_{i})$ be unit vectors such that $l_{1}=\mathbb{R}_+v_{1}$ and $l_{2}=\mathbb{R}_+v_{2}$. Obviously $l_{1}\neq l_{2} \Longleftrightarrow v_{1} \neq v_{2}$. Then we know that our curves are of the form $\gamma_{1}(t)=tv_{1}+g_{1}(t)$, $\gamma_{1}(t)=tv_{1}+g_{2}(t)$, where $||g_{i}(t)||=o(t)$.
Notice that
\[||\gamma_{1}(t)-\gamma_{2}(t)||=||(tv_{1}+g_{1}(t))-(tv_{2}+g_{2}(t))||=||t(v_{1}-v_{2})+g(t)||\]
where $g(t)=g_{1}(t)-g_{2}(t)$ and $||g(t)||=o(t)$.
We get that $tord(\gamma_{1},\gamma_{2}) \geq 1$ and $tord(\gamma_{1},\gamma_{2})=1 \Longleftrightarrow v_{1}\neq v_{2}$. Now consider the inner metric. First observe that $d_{l}(\gamma_{1}(t),\gamma_{2}(t)) \geq ||\gamma_{1}(t)-0||+||0-\gamma_{2}(t)|| \geq t+t=2t$. and since the pancake metric is equivalent to the length metric, we get that $tord_{inn}(\gamma_{1},\gamma_{2}) \leq 1$. We have two possible cases:\\
1) The germ is normally embedded. In this situation $d_{ind}(\gamma_{1}(t),\gamma_{2}(t))=d_{l}(\gamma_{1}(t),\gamma_{2}(t))$ so that $tord(\gamma_{1},\gamma_{2})=tord_{inn}(\gamma_{1},\gamma_{2})$ and from the considerations above we know that they both have to be equal to $1$ so that $v_{1}=v_{2}$.\\
2) The germ is not normally embedded i.e. it requires more than one pancake. Taking a neighbourhood of $0$ small enough we can assume that in the pancake decomposition we have only two pancakes. From the property 3) from Theorem \ref{pancake} we deduce that they can only intersect in $0$ and from the construction of the pancake metric we obtain that 
$d_{P}(\gamma_{1}(t),\gamma_{2}(t))=||\gamma_{1}(t)-0||+||0-\gamma_{2}(t)||=||\gamma_{1}(t)||+||\gamma_{2}(t)||=t+t=2t$ 
so $tord_{inn}(\gamma_{1},\gamma_{2})=1$. Then
\[(\gamma_{1} \cup \gamma_{2},0) \; \text{is normally embedded}\; \Longleftrightarrow tord(\gamma_{1},\gamma_{2})=1 \Longleftrightarrow v_{1}\neq v_{2},\]
which ends the proof.
\end{proof}
\end{theorem}

\begin{definition}
Let $X \subset \mathbb{R}^{n}$ be a semi-algebraic set. We define the dimension of $X$ as
$$ \dim X = \max \{\dim \Gamma\;|\; \Gamma \subset X,\; \Gamma - \text{a semi-algebraic submanifold of class}\; \mathcal{C}^{1} \}. $$
\end{definition}

\begin{definition}
For a semi-algebraic set $X \subset \mathbb{R}^{n}$ and point $a \in X$ we define the dimension of $X$ at $a$ as 
$$ \dim_{a} X = \min \{\dim X \cap U\;|\; U - \text{an open neighbourhood of}\ a \}. $$
Then the dimension of the germ $(X,0)$ is by definition $\dim_0 X$ as it does not depend on the choice of the representative.
\end{definition}

\begin{corollary}
\label{smooth}
A $\mathcal{C}^{1}$ semi-algebraic germ of dimension 1 in $\mathbb{R}^{2}$ is normally embedded.
\end{corollary}

\begin{remark}
The corollary above follows directly from the last theorem but in fact it is true for any definable smooth germ in $\mathbb{R}^{n}$. Indeed such a germ is a graph of a $\mathcal{C}^{1}$, hence Lipschitz function germ over its tangent space.
\end{remark}


 Before that, let us recall that given two semi-algebraic curve germs $\delta_i$ at $0\in\mathbb{R}^n$ that are distinct but have the same tangent half-line $l$, it follows from the general theory that we can choose a linear coordinate system so that $l$ is the non-negative $x$-axis and both $\delta_i$ are graphs over it with $\delta_1(t)<\delta_2(t)$ (after renumbering, if necessary) for all $t>0$ sufficiently small. In such a situation it becomes self explanatory what we mean when we say that a curve lies between $\delta_1$ and $\delta_2$ -- its germs lies between the two graphs.

\begin{proposition}
\label{side}
Suppose we have a one-dimensional semialgebraic set $X \subset \mathbb{R}^{2}$
and two semialgebraic curves $\delta_{1},\delta_{2}:[0,\epsilon) \rightarrow \mathbb{R}^{2}$, $\delta_{1}(0)=\delta_{2}(0)=0$, $(\delta_{1},0) \neq (\delta_{2},0)$ but sharing the same tangent half-line $l$ at the origin and $\delta_i\subset\overline{M_X}$, $i=1,2$. 

Then there exists a branch of $X$, i.e. a curve $\gamma\subset X$ through 0, that lies between $\delta_{1},\delta_{2}$, $\gamma(t) \in m(\delta_{i_{0}}(t))$ and so $C_0(\gamma)=l$.
\begin{proof}
 It is easy to see that for $a \in m(x)$, the vector $x-a$ is normal to $X$ at $a$ in the sense that for any $v\in C_{a}(X)$, the inner product $\langle x-a,v
\rangle\leq 0$. It follows that for any $z:=(x,y)$ between the two graphs, the closest point on $c:=(a,b)\in X$ is not the origin. We may restrict our considerations to a disc centred at the origin. Then if there were no point from $X$ between the two graphs, the segment joining $z$ and $c$ has to intersect one of the graphs $\delta_i$ at some point $w$. Then since $\#m(w)>1$ and $c\in m(w)$, we have a point $c'\in m(w)\setminus\{c\}$ and by the triangle inequality $||w-c'||<||w-c||$ which is a contradiction.

Once we know, by the argument above, that between the two graphs there must be a branch $\tilde{\gamma}$ of $X$, we necessarily have $C_0(\tilde{\gamma})=l$.
If $\tilde{\gamma}$ realizes the closest points to $\delta_{1}$ or $\delta_{2}$ then we are done. Otherwise, from Darboux property, between $\tilde{\gamma}$ and $\delta_{i}$ we would have another curve from the medial axis and we repeat the argument. By the semi-algebraicity such a situation can appear only finitely many times.
\end{proof}
\end{proposition}

\begin{proposition}
Let $(X,0) \subset \mathbb{R}^{2}$ be a germ of dimension 1 in $0$. If $(X,0)$ is not normally embedded, then $0 \in \overline{M_{X}} \cap X$.
\begin{proof}
The germ $(X,0)$ decomposes into finitely many curve branch germs $\{\gamma_{1},...,\gamma_{k}\}$. There must be $k \geq 2$ because otherwise we can extend the branch $\gamma_{1}$ to form a $\mathcal{C}^{1}$ curve. Then the Corollary \ref{smooth} shows that it is normally embedded. Therefore we have two branches $\gamma_{1},\gamma_{2}$ with no other branch in between (in the sense of the previous proposition) such that they share a common tangent half-line at 0 (Theorem \ref{curvegerm}).  Without loss of generality we can assume that the two curves are the graphs of two functions over the tangent half-line. The definability implies that we may assume these functions have constant convexity. 
As in the previous proof no point between $\gamma_{1},\gamma_{2}$ has its closest at $0 \in \gamma_{1} \cap \gamma_{2}$. Now taking the vertical sections between these two curves, by the Darboux property, we get a point that belongs to the medial axis and this arbitrarily near the origin.
\end{proof}
\end{proposition}

\begin{theorem}
\label{lastmedial}
Let $(X,0)$ be a closed semialgebraic germ of $\mathbb{R}^{2}$ and $M_{X}$ its medial axis. Assume that $0 \in \overline{M_{X}} \cap X$. Then
\[\overline{M_{X}}\; \text{is not normally embedded}\; \Longrightarrow X \; \text{is not normally embedded}.\]
\begin{proof}
From Theorem \ref{lne} we know there exists two parametrized by the distance semialgebraic curves $\delta_{1},\delta_{2} \subset \overline{M_{X}}$ such that $\delta_{1}(0)=\delta_{2}(0)=0$ and $tord(\delta_{1},\delta_{2})>tord_{inn}(\delta_{1},\delta_{2})$. Therefore, they share a common tangent half-line $l$. Then for every $t>0$ we have $\# m(\delta_{1}(t))>1$ so that there exists $x_{t}^{(1)},x_{t}^{(2)} \in m(\delta_{1}(t)), \; x_{t}^{(1)} \neq x_{t}^{(2)}$. Thanks to the Curve Selecting Lemma we obtain two semialgebraic curves $\gamma_{1},\gamma_{2} \subset X$ such that $\gamma_{i}(t)=x_{t}^{(i)}$; denote by $\tilde{\gamma_{i}}$ their reparametrizations by the distance.\\
By Proposition \ref{side} we may assume that at least one of these curves has also the same tangent half-line $l$. Now, since we know that $\delta_{1}$ and $\gamma_{1}$ have a common tangent half-line, arguing similarly as in the proof of the previous theorem, we get that $tord(\delta_{1},\tilde{\gamma_{2}})>1$. Since $||\delta_{1}(t)-\tilde{\gamma_{1}}(t)||\leq||\delta_{1}(t)-\gamma_{1}(t)||$ we have that $tord(\delta_{1},\gamma_{1})>1$.\\
Notice that $||\delta_{1}(t)-\gamma_{1}(t)||=||\delta_{1}(t)-\gamma_{2}(t)||$, so that $tord(\delta_{1},\gamma_{1})=tord(\delta_{1},\gamma_{2})$. This means that $tord(\delta_{1},\gamma_{2})>1$. Let $v \in C_{0}(\delta_{1})$ be a unit vector such that $l=\mathbb{R}_+v$. Now
\[\frac{||tv-\gamma_{2}(t)||}{t}\leq \frac{||tv-\delta_{1}(t)||+||\delta_{1}(t)-\gamma_{2}(t)||}{t}=\frac{||tv-\delta_{1}(t)||}{t}+\frac{||\delta_{1}(t)-\gamma_{2}(t)||}{t} \rightarrow 0,\]
so we get that $\gamma_{2}$ is also tangent to the line $l$. This means that $\gamma_{1}$ and $\gamma_{2}$ have a common tangent since tangents to curves are unique. From Theorem \ref{curvegerm} we get that the germ $(\gamma_{1} \cup \gamma_{2},0)$ is not normally embedded so, based on Theorem \ref{lne}, the germ $(X,0)$ cannot be normally embedded (we have only to notice that since $\gamma_{1}$ is `separated' from the other parts of the set, the pancake metric for $\gamma_{1}$, $\gamma_{2}$ considered as a germ of two curves is the same as if we consider them as a part of $X$ -- in such a case the pancake metric possibly could be smaller but is not because of the behaviour of the medial axis).  
\end{proof}
\end{theorem}

\begin{remark}
From the proof of the last theorem we get that for $(X,0) \subset \mathbb{R}^{2}$, we have
\[\mathcal{L}(\overline{M_{X}},0) \leq \mathcal{L}(X,0).\]
\end{remark}

\end{document}